\documentclass[10pt,a4paper]{elsarticle}

\usepackage{atbegshi}
\AtBeginDocument{\AtBeginShipoutNext{\AtBeginShipoutDiscard}}

\usepackage{proof}
\usepackage{amsthm}

\newtheorem{theorem}{Theorem}
\def\ie{{ i.e. }}

\usepackage[latin2]{inputenc}
\usepackage{amsmath}
\usepackage{amsfonts}
\usepackage{amssymb}
\usepackage{makeidx}

\def\nfrac#1#2{\mbox{\footnotesize $\displaystyle\frac{#1}{#2}$}}

\journal{arXiv.org}

\begin{document}

\setcounter{page}{0}

\begin{frontmatter}

\title{\bf \large A PROOF OF AN OPEN PROBLEM OF YUSUKE NISHIZAWA FOR A POWER-EXPONENTIAL FUNCTION}

\author[etf]{\bf Branko Male\v sevi\' c\corref{cor}}
\ead{branko.malesevic@etf.rs}
\author[etf]{\bf Tatjana Lutovac}
\ead{tatjana.lutovac@etf.rs}
\author[etf,ftn]{\bf Bojan Banjac},
\ead{bojan.banjac@uns.ac.rs}

\cortext[cor]{Corresponding author.                                                 \\
Authors were supported in part by Serbian Ministry of Education,
Science and Technological Development, Projects ON 174032, III 44006 and TR 32023.}

\address[etf]{Faculty of Electrical Engineering, University of  Belgrade,           \\
Bulevar kralja Aleksandra 73, 11000 Belgrade, Serbia}

\medskip

\address[ftn]{University of Novi Sad, Faculty of Technical Sciences,                \\
Computer Graphics Chair, Trg Dositeja Obradovi\' ca 6,                              \\
21000 Novi Sad, Serbia}

\vspace*{-5.0 mm}
\begin{abstract}
{\small This paper presents a proof of the following conjecture, stated by Yusuke Nishi\-zawa in [Appl. Math. Comput. 269, (2015), 146--154.]:
for $\displaystyle 0<x<\pi/2$ the inequality \mbox{\small $\displaystyle \frac{\sin{x}}{x} \!>\! \left(\frac{2}{\pi} +
\frac{\pi\!-\!2}{\pi^{3}}(\pi^{2}\!-\!4x^{2})\right)^{\!\theta(x)}\!$} holds, where \mbox{\small $ \displaystyle \theta(x) \! = \!
-\frac{(48\!-\!24\pi\!+\!\pi^{3})x^{3} }{3(\pi\!-\!2)\pi^{3}}+\frac{\pi^{3}}{24(\pi\!-\!2)}. $}}
\end{abstract}


\end{frontmatter}

\section{Introduction}

In \cite{Nishizawa_2015}, Nishizawa  proved the following power-exponential~inequalities:

\begin{theorem}
$\,$ For  $\displaystyle 0<x<\pi/2$, we have

$\displaystyle \left( \frac{2}{\pi} +
\frac{1}{\pi^{3}}(\pi^{2}-4x^{2}) \right)^{\theta}  <
\frac{\sin{x}}{x} < \left( \frac{2}{\pi} +
\frac{1}{\pi^{3}}(\pi^{2}-4x^{2}) \right)^{\vartheta} $

\noindent
with the best possible constants $\theta = 1$ and $\vartheta = 0$.
\end{theorem}
\begin{theorem} $\,$ For  $\displaystyle 0<x<\pi/2$, we have

$\displaystyle \left( \frac{2}{\pi} + \frac{\pi -2}{\pi^{3}}
(\pi^{2}-4x^{2}) \right)^{\theta}  < \frac{\sin{x}}{x} < \left(
\frac{2}{\pi} + \frac{\pi -2}{\pi^{3}}(\pi^{2}-4x^{2})
\right)^{\vartheta} $

\noindent
with the best possible constants $\theta =
\pi^{3}/(24(\pi-2))\cong1.13169$ and $\vartheta = 1$.
\end{theorem}
\begin{theorem} $\,$ For  $\displaystyle 0<x<\pi/2$, we have

$\displaystyle  \frac{\sin{x}}{x} < \left( \frac{2}{\pi} +
\frac{1}{\pi^{3}}(\pi^{2}-4x^{2}) \right)^{\theta(x)} $

\noindent
where $\theta(x) = 4x^{2}/\pi^{2}$.
\end{theorem}

{\bf The open problem of Yusuke Nishizawa.}
Considering the previous theorems, Nishizawa stated  the following
open problem~(Problem ~3.1~of~\cite{Nishizawa_2015}):

\break

\noindent
{\em For $\, \displaystyle 0<x<\pi/2$, we have

\begin{equation}
\label{startinequality0}
 \displaystyle \frac{\sin{x}}{x}
> \left( \frac{2}{\pi} + \frac{\pi -2}{\pi^{3}}(\pi^{2}-4x^{2})
\right)^{\theta(x)}
\end{equation}

\noindent
where  $\theta(x)$ is the function of $x$ and  $\, \displaystyle
\theta(x)= -\frac{(48-24\pi+\pi^{3})x^{3} }{3(\pi-2)\pi^{3}}+
\frac{\pi^{3}}{24(\pi-2)}$.}

\bigskip
This paper provides a proof of  Nishizawa's open problem, using
approximations and methods  from \cite{Mortici_2011} and
\cite{Malesevic_Makragic_2015}. Let us notice that method from paper \cite{Malesevic_Makragic_2015} relates
to proving mixed trigonometric polynomial inequalities
\begin{equation}
\label{MTF 1}
f(x) = \sum\limits_{i=1}^{n}\alpha_{i}x^{p_{i}}\!\cos^{q_{i}}\!x\sin^{r_{i}}\!x > 0,
\end{equation}
where $\alpha_{i} \!\in\! \mathbb{R}$, \mbox{$p_{i}, q_{i}, r_{i} \!\in\!
\mathbb{N}_{0}$}, \mbox{$n \!\in\! \mathbb{N}$}~and~$x \!\in\! (\delta_{1},\delta_{2})$, $\delta_1 \!\leq\! 0
\!\leq\!\delta_2$, with $\delta_{1} \!<\! \delta_{2}$. The function $f(x)$ is called mixed trigonometric polynomial function.

\smallskip
In this paper the power-exponential inequality (\ref{startinequality0}) can be rewritten as one example of the inequality
of the following form:
\begin{equation}
\label{Funct_F}
F(x)
=
f(x)
+
\sum\limits_{j=1}^{m}P_{j}(x)\ln\!\left(f_{j}(x)\right) > 0,
\end{equation}
where $f(x)$ and $f_{j}(x)$ are mixed trigonometric polynomial functions, with $f_{j}(x)\!>\!0$
for $x \in \left(0, \frac{\pi}{2}\right)$; $P_{j}(x)$ is a real po\-ly\-nomial of degree $k_{j}$ and
\mbox{$m \!\in\! \mathbb{N}$}.$\;$ The function $F(x)$ is called mixed logarithmic-trigonometric polynomial
function and the inequality $F(x) > 0$ is called mixed logarithmic-trigonometric polynomial inequality.

\smallskip
Assuming that in mixed logarithmic-trigonometric polynomial function appear only mo\-no\-mials and logarithmic functions,
and no trigonometric sine and cosine functions, we can call that function mixed logarithmic polynomial function.
One application of mixed logarithmic polynomial functions, in purpose of proving some power-exponential inequalities,
was given in the papers~\cite{Nishizawa_2015},~\cite{Nishizawa_2015b},~\cite{Nishizawa_2015c}.

\smallskip
Let us assume that the degree of the zero polynomial is $-1$.
Then, for a mixed logarithmic-trigonometric polynomial function $F(x)$ it is not difficult to show the following result:
\begin{equation}
\label{Properties}
\!\!
\begin{array}{l}
\mbox{\em The$\;\;$derivative$\;\;F^{(\mbox{\tiny $K$}\!+\!1)\!\,}(x)\;\;$is$\;\;$quotient$\;\;$of$\;\;$two$\;\;$mixed$\;\;$trigono-} \\
\mbox{\em metric$\;\;$polynomial$\;\;$functions,$\;\,$where$\;\,\mbox{\small $K$}\!=\!\max\{k_{j}\,|\,j\!=\!1,\!...,m\}$.}
\end{array}
\end{equation}

\smallskip
In the proof of Nishazawa's open problem, we will use the following statement.

\begin{theorem}
\label{TH}
Let $\Phi\!:\!\left(0, \delta\right) \!\longrightarrow\! \mathbb{R}$ $\left(\delta \!>\!0\right)$,
be a $(\mbox{\small $K$}\!+\!1)$-times differentiable function $\left(\mbox{\small $K$} \!\in\!
\mathbb{N} \cup \{0\}\right)$ such that

\medskip
\noindent
$\,\;(i)\,\,$ $\displaystyle\lim_{x \rightarrow 0+}{\Phi(x)} \geq 0$,

\smallskip
\noindent
$\,(ii)\,$ $\displaystyle\lim_{x \rightarrow 0+}{\Phi^{(j)}(x)} \geq 0$ holds true for every $j \in \{1,2,\ldots,\mbox{\small $K$}\}$,

\smallskip
\noindent
$(iii)$ $\Phi^{(\mbox{\tiny $K$}+1)}(x) > 0$ for $0 \!<\!x\!<\! c$ where $c \in \left(0,\delta\right]$.

\medskip
\noindent Then, for $0 \!<\!x\!<\! c$, inequality
\begin{equation}
\label{Ineq_2}
\Phi(x)
>
0
\end{equation}
holds true.
\end{theorem}

\break

\smallskip
In the next section the proof of Nishizawa's open problem also makes use of the
fact that for the constant $\pi$ and  a given rational function $R(x)$, it is
possible to determine either $R(\pi) > 0$ or $R(\pi) < 0$. Stated is a consequence
of the fact that for an arbitrarily small $\varepsilon > 0$,  there exist fractions
$p/q$ and $r/s$ such that $p/q > \pi > r/s$ and $p/q - r/s < \varepsilon$. Fractions
$p/q$ and $r/s$ can be chosen as two consequential convergents in the development of
continued fractions of $\pi$.

\section{The proof of Nishizawa's open problem}

As $\, \displaystyle \nfrac{\sin{x}}{x} > 0$ and $\,\,\displaystyle \nfrac{2}{\pi} + \nfrac{\pi -2}{\pi^{3}}(\pi^{2}
-4x^{2}) > 0 \, \,$ for $ x \in (0,\pi/2)$, the power-exponential inequality (\ref{startinequality0}) is equivalent to the following
inequality:

\begin{equation}
\label{startinequality02}
 \displaystyle \ln \frac{\sin{x}}{x}  \,  >  \,
 \ln \left( \frac{2}{\pi} + \frac{\pi -2}{\pi^{3}}(\pi^{2}-4x^{2})
\right)^{\theta(x)}
\end{equation}

\noindent or

\noindent
\begin{equation}
\label{startinequality03}
\displaystyle \ln\frac{\sin{x}}{x}  \, -  \,
\theta(x) \ln \left( \frac{2}{\pi} + \frac{\pi
-2}{\pi^{3}}(\pi^{2}-4x^{2}) \right)  \, > \, 0,
\end{equation}

\noindent
where $\displaystyle \theta(x)= -\frac{(\pi^{3}-24\pi+48)x^{3} }{3(\pi-2)\pi^{3}} + \frac{\pi^{3}}{24(\pi-2)}$.

\bigskip
\noindent
Let us notice that the previous inequality is a mixed logarithmic-trigonometric polynomial inequality

\begin{equation}
F(x)
=
\ln \sin x - \ln x - \theta(x) \ln \! \left( \frac{2}{\pi} + \frac{\pi\!-\!2}{\pi^{3}}(\pi^{2}\!-\!4x^{2})\right) \!>\!0
\end{equation}
for $x \in  \left(0, \pi/2 \right)$.

\medskip
Let us further consider the following mixed logarithmic-trigonometric polynomial function

\begin{equation}
F_{1}(x)
=
\displaystyle \ln \sin x - \ln x  \, -  \,
\theta_{1}(x) \ln
\! \left( \frac{2}{\pi} + \frac{\pi-2}{\pi^{3}}(\pi^{2}-4x^{2}) \right)
\end{equation}

\medskip
\noindent
where $\,\displaystyle \theta_{1}(x)=\frac{(\pi^{3}-60\pi+120)}
{720(\pi-2)}x^{2} + \frac{\pi^{3}}{24(\pi-2)}$.

\medskip
\noindent
As

\begin{equation}
\label{ineq_(10)}
\frac{2}{\pi} + \frac{\pi-2}{\pi^{3}}(\pi^{2}-4x^{2})
\, < \,
\frac{2}{\pi} + \frac{\pi -2}{\pi^{3}}\pi^{2}
\, = \,
1
\end{equation}

\medskip \noindent
we can conclude that for $x \in  \left(0, \pi/2 \right)$:
\begin{equation}
\begin{array}{cl}
                     &  F(x)- F_{1}(x) \, \geq 0                                              \\[3.0 ex]
 \Longleftrightarrow & \theta(x) - \theta_{1}(x) \geq 0                                       \\[2.0 ex]
 \Longleftrightarrow \, &  \displaystyle - \, \frac{1}{3(\pi -2)}
 \left(\frac{\pi^{3}-60\pi+120}{240} + \frac{\pi^{3}-24\pi+48}{\pi^{3}}x\right) x^{2} \geq 0. \\
\end{array}
\end{equation}

\noindent
As the last inequality holds for $\,  x \in (0, c]$, ~ where ~ $\displaystyle c
=
-\frac{ \left( {\pi }^{3}-60(\pi-2) \right) {\pi }^{3}}{240 \,({ \pi}^{3}-24\,(\pi-2))}= 1.342\ldots$
we distinguish two cases $x \in (0, c]$ or $x \in (c, \pi/2)$.

\break

\subsection{The case 1: $x \in \left(0, c\,\right]$}

\medskip

It is enough to prove that
$
F_{1}(x) > 0
$
for $x\in(0,c]$. We have:

\begin{equation}
\label{inequality}
F_{1}'''(x)
\, = \,
\frac{~ 2 \, A(x) \, \sin^{3}{x} \,\, + \,\, B(x)  \cos{x}~}{ 45 \, x^{3} \, C(x) \, \sin^{3}{x}}
\end{equation}
where
\begin{equation}
\begin{array}{rcl}
C(x) \!\!&\!\!=\!\!&\!\!  -64\,(\pi -2)^{3} {x}^{6} \,   +  \, 48 {\pi}^{3}(\pi -2)^{2}{x}^{4} \, -  \, 12 \pi^{6}(\pi-2)x^{2} \, + \pi^{9}
                                                                                \\
     \!\!&\!\! \!\!&\!\!                                                        \\
     \!\!&\!\!=\!\!&\!\! (\pi^{3}- 4 (\pi-2)x^{2} )^{3},
\end{array}
\end{equation}

\medskip

\begin{equation}
B(x) \,=\, \displaystyle 90 x^{3}\left(\pi^{3}\, - \, 4(\pi
-2)x^{2}\right)^{3} \,= \, 90 \, x^{3} \, C(x) ~~~~~
\end{equation}

\noindent and

\begin{equation}
\begin{array}{rcl}
A(x)
\!\!&\!\!=\!\!&\!\! 8(\pi-2)^{2}(\pi^{3}-60\pi+120)x^{8}                                        \\[2.0 ex]
\!\!&\!\! \!\!&\!\! -6 (\pi-2)(\pi^{6}-100\pi^{4}+200\pi^{3}-480(\pi-2)^{2})x^{6}               \\[2.0 ex]
\!\!&\!\! \!\!&\!\! + 3\pi^{3} (\,\pi^{6}-720 (\pi -2)^{2} \,)x^{4}                             \\[2.0 ex]
\!\!&\!\! \!\!&\!\! + 540 \pi^{6}(\pi-2)x^{2}                                                   \\[2.0 ex]
\!\!&\!\! \!\!&\!\! -45 \pi^{9}.
\end{array}
\end{equation}

\medskip
\noindent
Let us determine the sign of the polynomials $\, C(x),\, B(x)$ and $\, A(x)$ for $x \in (0, c]$.\\
By substituting $\, t= 4 (\pi-2)x^{2}\,$   for $x \in  (0, c]$, the polynomial $\,C(x)$  can
be transformed into  the  polynomial $\displaystyle ~ C_{1}(t) \, = \, (\pi^{3}- t)^{3}$
for $\,\, t \in (0,\, 4(\pi-2)c^{2}]$. Obviously, the sign of the polynomial $C_{1}(t)$
coincides with the sign of the polynomial $\displaystyle \pi^{3}-t$ for $\, t \in (0,\,
4(\pi-2)c^{2}]$.

\medskip
\noindent
Since  $\, \displaystyle \left(\frac{\pi^{3}-60(\pi-2)}{\pi^{3}-24(\pi-2)}\right)^{2} \, < \, 1$,
we have
\begin{equation}
\begin{array}{rcl}
\displaystyle \pi^{3}- 4(\pi-2)c^{2}
 &\!\!=\!\!&
\displaystyle  \, \pi^{3} \, - \,
4(\pi-2)\left(\frac{(\pi^{3}-60(\pi-2))}{240(\pi^{3}-24(\pi-2))}
\pi^{3} \right)^{2} \,                                                                      \\[-1.00 ex]
 &\!\! \!\!&                                                                                \\[-1.00 ex]
 &\!\!=\!\!&  \pi^{3} \displaystyle \left( 1 -
(\pi-2)\frac{1}{14400}\left(\frac{(\pi^{3}-60(\pi-2))}{(\pi^{3}-24(\pi-2))}\right)^{2}
\pi^{3}  \right)                                                                            \\[-1.00 ex]
 &\!\! \!\!&                                                                                \\[-1.00 ex]
 &\!\!>\!\!& \pi^{3} \displaystyle \left(  \displaystyle 1 -
(\pi-2)\frac{1}{14400}
\pi^{3} \,\right)  \,                                                                       \\[-1.00 ex]
 &\!\! \!\!&                                                                                \\[-1.00 ex]
 &\!\!=\!\!& \displaystyle \pi^{3} \left( \frac{14400 - (\pi-2)\pi^{3}}{14400}\right) \,\,  \\[-1.00 ex]
 &\!\! \!\!&                                                                                \\[-1.00 ex]
 &\!\!>\!\!& \, 0.
\end{array}
\end{equation}

\noindent
Therefore, we can conclude  that $ C_{1}(t)>0$ for  $t \in
( \, 0, \, 4(\pi-2)c^{2}\, ]
 \subset (0, \pi^{3})$,   \ie
$C(x) \, > \, 0$  for  $x \in (0, \, c]$  and $B(x) \,
> \, 0$ for  $x \in (0, \, c]$.

\smallskip\noindent
Let us prove that
\begin{equation}
\label{A}
A(x)<0,
\end{equation}
for $x\in (0,\, c]$. We note  that $A(x)$ can be written as
\begin{equation}
A(x)\,= \,2(\pi-2)x^{6}\varphi_{1}(x) \,\, + \,\, 3\pi^{3} x^{2}
\varphi_{2}(x) \,\, - \,\,45 \pi^{9},
\end{equation}
where
\begin{equation}
\displaystyle \begin{array}{rcl}
\varphi_{1}(x)  &\!\!\!\!\!=4(\pi\!-\!2)(\pi^{3}\!-\!60\pi\!+\!120)x^{2}
\!-\!
3\left(\pi^{6}\!-\!100\pi^{4}\!+\!200\pi^{3}\!-\!480(\pi\!-\!2)^{2}\right)   \\
\end{array}
\end{equation}
and
\begin{equation}
\begin{array}{rcl}
\varphi_{2}(x) \!&\!\!=\!\!&\!(\pi^{6}\!-\!720(\pi\!-\!2)^{2}\,) x^{2}\,+\,180\pi^{3}(\pi\!-\!2).    \\
\end{array}
\end{equation}

\noindent
As $\displaystyle \, \pi^{3}-60\pi-120 <0$ and
$\displaystyle \,\pi^{6}-720(\pi-2)^{2} \,> \,0$, we have the
following estimation, for $\, x\in (0, c]$:
\begin{equation}
A(x) \, \leq \, \,2(\pi-2)c^{6}\varphi_{1}(0) \,\, + \,\, 3\pi^{3}
c^{2} \varphi_{2}(c) \,\, - \,\,45 \pi^{9} \, =
-138097.868\ldots \, < \, 0.\end{equation}

\noindent
In view  of all the above, we can conclude that for $x \in (0, \, c]$:
\begin{equation}
C(x) > 0, ~~ B(x)>0 ~~ \mbox{and}
~~ A(x)<0.
\end{equation}
Now we prove that
\begin{equation}
\label{equationx}
g_{1}(x)
=
2\, A(x)\, \sin^{3}{x} + B(x) \,\cos{x} \, > \, 0
\end{equation}
for $x\in(0,c]$. Let us note that $g_{1}(x)$ is a mixed trigonometric polynomial function,
and that the proof of previous inequality will be proved applying the methods from
\cite{Mortici_2011} and \cite{Malesevic_Makragic_2015}. In particular,  we use the following
inequalities from \cite{Malesevic_Makragic_2015}:
\begin{equation}
\cos x > 1-\frac{x^{2}}{2!}\, + \, \frac{x^{4}}{4!}\, - {\frac {{x}^{6}}{6!}},
\qquad {\big (}x \in (0,c)\subseteq (0,\sqrt{90}){\big )};
\end{equation}
and
\begin{equation}
\sin x >  x-\frac{x^{3}}{2!}+{\frac {{x}^{5}}{5!}}  , \qquad
{\big (}x \in (0,c)\subseteq (0,\sqrt{72}){\big )}.
\end{equation}
Therefore,  for $x \in \left(0, \pi/2\right)$,  we have:
\begin{equation}
\begin{array}{rcl}
g_{1}(x) \,=\, 2 A(x) \sin^{3}{x} + B(x) \cos{ x}
&\!\!>\!\!&
2 A\left(x\right) \left( x-\nfrac{x^{3}}{2!} \!+\! {\nfrac{{x}^{5}}{5!}} \right)^{3}              \\[1.5 ex]
&\!\! \!\!& \!\!+\,
B\left(x\right)  \left( 1\!-\!\nfrac{x^{2}}{2!} \!+\! \nfrac{x^{4}}{4!} \!-\! {\nfrac {{x}^{6}}{6!}} \right)
                                                                                                  \\[1.0 ex]
&\!\!=\!\!& \nfrac{x^{9}}{864000} \, P_{14}(x)\,,
\end{array}
\end{equation}

\noindent
where $P_{14}(x)$ is the polynomial of the $14^{th}$ degree, as follows:

\vspace*{-2.0 mm}

{\small
\begin{equation}
\begin{array}{rl}
\! P_{14}(x) \!\! & \! \! =  \!   \displaystyle \sum_{i=0}^{14} a_{i} x^{i}        \\[1.50 ex]
& = 8\! \left( {\pi }^{3}-60\!\pi +120 \right) \left( \pi -2 \right)^{2} x^{14}    \\[1.50 ex]
& - \! 6\! \left( \pi -2 \right) \left( {\pi }^{6}-20{\pi}^{4}+40{\pi }^{3}-5280{\pi }^{2}+21120\pi -21120 \right)x^{12}
                                                                                   \\[1.50 ex]
& + \! 3\!\left( {\pi }^{9}+ 40\! \left( \pi -2 \right)
\left(3{\pi }^{6}-214{ \pi }^{4}+428{\pi }^{3}-7680{\pi}^{2}+30720\pi -30720 \right) \right) x^{10}
                                                                                   \\[1.50 ex]
& - \! 20 \! \left(\! 9\pi^{9}\! + \!(\pi-2)\left(  441\pi^{6}\! - \!44320\pi^{4}\! + \! 88640\pi^{3} \! - \!762240\pi^{2} \!+ \! 3048960\pi  \! - \! 3048960 \right)\! \right) \! x^{8}
                                                                                   \\[1.50 ex]
& + \! 15 \! \left( \! 309 \pi^{9} \! + \! 400 \! (\pi -2) \! \left( 17\pi^{6} \! - \! 2552 \pi^{4}  \! + \! 5104 \pi^{3}  \! - \! 24576 \pi^{2} \! + \! 98304 \pi - 98304  \right)     \! \right) \! x^{6}
                                                                                   \\[1.50 ex]
& -  \! 300  \! \left( \! 215\pi^{9} \! + \! 72 (\pi-2) \left(13\pi^{6} - 6880 \pi^{4} + 13760 \pi^{3} - 34560 \pi^{2} + 138240\pi - 138240  \right)   \! \right) x^{4}
                                                                                   \\[1.50 ex]
& + 5400 \! \left( 91\pi^{9} \! - \! 80 (\pi-2)\left( 13\pi^{6} + 1744\pi^{4}-3488\pi^{3} + 1920\pi^{2} -7680\pi + 7680 \right) \! \right) x^{2}
                                                                                   \\[1.50 ex]
& -36000{\pi }^{3} \! \left( 47\!{\pi }^{6}-1440\! \left(\pi -2\right)\left( {\pi }^{3}+20\!\pi -40 \right)  \right).
\end{array}
\end{equation}
}

\noindent Therefore, for  inequality (\ref{equationx}) it is
sufficient to prove that
\begin{equation}
\label{P14}
\displaystyle\,P_{14}(x)\,> \,0,
\end{equation}
for $x\in (0,c]$. It is easy to
check that non-zero coefficients $a_{i}, \,
i\in\{14, 12, 10, 8, 6,$ $4, 2, 0\} $ of the polynomial $P_{14}(x)$ satisfy
the following conditions:  $a_{14}<0$, \mbox{$a_{12}>0$},  $a_{10}<0$, $a_{8}>0$, $a_{6}<0$,
$a_{4}>0$, $a_{2}<0$ and $a_{0}>0$.  Thus, for the proof of  $\,\, P_{14}(x)\, = \,\,
x^{12} \left(a_{14} x^{2} + a_{12}\right) \, \, +  x^{8}
\left(a_{10} x^{2} + a_{8}\right)\, \, +  x^{4} \left(a_{6} x^{2} +
a_{4}\right)$

\noindent $ \, \, +
  \left(a_{2} x^{2} + a_{0}\right) \, \,
>\,0$,  for $\,x \in (0,c]$,  it is sufficient and easy  to check
that the following inequalities  hold:  $\,a_{14}\, c^{2} \, + \,
a_{12} \,
> \, 0, ~~~ a_{10} c^{2} + a_{8} \, > \, 0, ~~~ a_{6} c^{2} + a_{4} \,
> \, 0$ and $\,a_{2} \,c^{2}\, + \, a_{0} \, > \, 0$.
Thus we may conclude that $\,P_{14}(x)>0$~and~based~on $F_{1}'''(x) \!=\! g_{1}(x) /
\!\left(45 \, x^{3} \, C(x) \, \sin^{3}{x}\right)
\!>\! P_{14}(x) \, x^{9} / \left(38880000  \, x^{3} \, C(x) \, \sin^{3}{x}\right)$
the inequality
\begin{equation}
F_{1}'''(x) > 0
\end{equation}
was proved for $\,x \in (0, c]$. Let us notice that
\begin{equation}
\displaystyle \lim_{x\rightarrow +0} F_{1}''(x)= 0,
\end{equation}
\begin{equation}
\displaystyle \lim_{x\rightarrow +0} F_{1}'(x)= 0
\end{equation}
and
\begin{equation}
\displaystyle \lim_{x\rightarrow +0} F_{1}(x)= 0.
\end{equation}
Then based on Theorem \ref{TH} we can conclude that
\begin{equation}
F_{1}(x) > 0 ~~ \mbox{for} ~~ x \in (0, c]
\end{equation}

\smallskip
\noindent
which also proves that $\,F(x) \,> \,0$  for $x\in (0,c]$.

\smallskip

\subsection{The case 2: $x \in \left(c, \pi/2\right)$}

\bigskip

In this subsection we prove that
$
F(x) >  0
$
for $x\in\left(c,\nfrac{\pi}{2}\right)$.
Let us consider the mixed logarithmic-trigonometric polynomial function
\begin{equation}
\!\!\!
G(x) \!=\! F(\nfrac{\pi}{2}\mbox{\small $-$}x) \!=\! \displaystyle \ln \cos x \mbox{\small $-$} \ln \!\left(\nfrac{\pi}{2} \!-\! x\right)
\!\mbox{\small $-$}\,
\omega(x) \ln \! \left(\nfrac{2}{\pi}
\mbox{\small $+$}
\nfrac{\pi\mbox{\small $-$}2}{\pi^{3}}\left(\pi^{2}\mbox{\small $-$}4\left(\nfrac{\pi}{2} \mbox{\small $-$} x\right)^{\!\!2}\right)\!\right)\!\!,
\end{equation}
where $\omega(x) = \theta{\big (}\frac{\pi}{2}-x{\big )} = -
{\nfrac {\left({\pi }^{3}-24\,\pi+48 \right)  \left(
\frac{\pi}{2}-x \right)^{3}}{ 3\left(\pi-2 \right) {\pi
}^{3}}}+{\nfrac{{\pi }^{3}}{24\,\pi-48}}$ and $x \in
\left(0,\nfrac{\pi}{2}\right)$. We have to prove the following mixed logarithmic-trigonometric polynomial inequality
\begin{equation}
\label{G(x)>0}
G(x) > 0
\end{equation}
for $x \in \left(0, c_1 \right)$, where \mbox{$c_1 = \nfrac{\pi}{2}
- c = 0.228 \ldots\;$}. Let us further consider the new mixed logarithmic-trigonometric polynomial function
\begin{equation}
G_{1}(x) = \displaystyle \ln \cos x - \ln\!\left(\frac{\pi}{2}\!-\!x\right)
-
\omega_{1}(x) \ln \! \left( \frac{2}{\pi}
+ \frac{\pi\!-\!2}{\pi^{3}}\left(\pi^{2}\!-\!4\left(\frac{\pi}{2}\!-\!x\right)^{2}\right) \right)\!,
\end{equation}
where $\omega_{1}(x) = \nfrac{x}{5}+1$ and $x \!\in\! \left(0,c_{1}\right)$.
Based on the inequality (\ref{ineq_(10)}) we can conclude that for $x \!\in\! \left(0,c_{1}\right)$:
\begin{equation}
\begin{array}{rl}
                    & G(x) > G_{1}(x)                                                             \\[2.0 ex]
\Longleftrightarrow & \omega(x) > \omega_{1}(x)                                                   \\[1.5 ex]
\Longleftrightarrow & \displaystyle\frac{x}{ \left( \pi -2 \right)
{\pi }^{3}}\,
{\Big (} \left( 20\,{\pi }^{3}-480\,\pi +960 \right) {x}^{2}                                      \\[1.5 ex]
                    &
\;\;\;\;\;\;\;\;\;\;\;\;\;\;\;\;\;
+\left( -30\,{\pi }^{4}+720\,{\pi }^{2}-1440\,\pi \right) x                                       \\[1.5 ex]
                    &
\;\;\;\;\;\;\;\;\;\;\;\;\;\;\;\;\; +\;15\,{\pi }^{5}-12\,{\pi
}^{4}-336\,{\pi }^{3}+720\,{\pi }^{2}\,{\Big )}> 0.
\end{array}
\end{equation}

\medskip
\noindent
As the last inequality holds for $x \!\in\! \left(0,c_{1}\right)$, in order to prove (\ref{G(x)>0}), it is enough to prove
\begin{equation}
G_1(x)>0,
\end{equation}
for $x \in (0,c_1)$. We have:
\begin{equation}
G_{1}''(x)=\displaystyle\frac{P(x)\cos^2\!x - \sin^2\!x \, Q(x)}{Q(x)\cos^2\!x}
\end{equation}
where
\begin{equation}
\begin{array}{rcl}
P(x) \!\!&\!\!=\!\!&\!\!
\left(-80\,{\pi }^{2}+320\,\pi -320\right) {x}^{6}                                                            \\[1.25 ex]
\!\!&\!\!+\!\!&\!\!
\left( 240\,{\pi }^{3}-992\,{\pi }^{2}+1088\,\pi -128 \right) {x}^{5}                                         \\[1.25 ex]
\!\!&\!\!+\!\!&\!\!
\left( -260\,{\pi }^{4}+1216\,{\pi }^{3}-1344\,{\pi }^{2}-576\,\pi +960 \right) {x}^{4}                       \\[1.25 ex]
\!\!&\!\!+\!\!&\!\!
\left( 120\,{\pi }^{5}-728\,{\pi }^{4}+720\,{\pi }^{3}+1472\,{\pi}^{2}-1920\,\pi \right) {x}^{3}              \\[1.25 ex]
\!\!&\!\!+\!\!&\!\!
\left( -20\,{\pi }^{6}+212\,{\pi }^{5}-132\,{\pi }^{4}-1184\,{\pi }^{3}+1440\,{\pi }^{2} \right) {x}^{2}      \\[1.25 ex]
\!\!&\!\!+\!\!&\!\!
\left( -24\,{\pi }^{6}-16\,{\pi }^{5}+408\,{\pi }^{4}-480\,{\pi }^{3} \right) {x}                             \\[1.25 ex]
\!\!&\!\!+\!\!&\!\! 11\,{\pi }^{6}-52\,{\pi }^{5}+60\,{\pi }^{4}
\end{array}
\end{equation}
and
\begin{equation}
Q(x)=5\, \left( \pi -2\,x \right) ^{2} \left(  \left( -2\,\pi +4
\right) {x }^{2}+ \left( 2\,{\pi }^{2}-4\,\pi  \right) x+{\pi }^{2}
\right) ^{2}.
\end{equation}
Obviously, $Q(x)>0$ for $x \in (0,c_1)$. Let us prove that $P(x)>0$
for $x \!\in\! (0,c_1)$. Note that $P(x)$ can be written as
\begin{equation}
P(x)
=
\phi_1 \left( x \right) +4\,{x}^{3}
\left( \pi -2 \right)  \left(  \left( -20\,\pi +40 \right) {x}^{3}+\phi_2 \left( x \right)  \right),
\end{equation}
where
\begin{equation}
\begin{array}{rcl}
\phi_1(x) \!\!&\!\!=\!\!&\!\!
\left( -20\,{\pi }^{6}+212\,{\pi }^{5}-132\,{\pi }^{4}-1184\,{\pi }^{3}+1440\,{\pi }^{2} \right) {x}^{2}      \\[1.25 ex]
\!\!&\!\!+\!\!&\!\!
\left( -24\,{\pi }^{6}-16\,{\pi }^{5}+408\,{\pi }^{4}-480\,{\pi }^{3} \right) {x}                             \\[1.25 ex]
\!\!&\!\!+\!\!&\!\! \left( 11\,{\pi }^{6}-52\,{\pi }^{5}+60\,{\pi
}^{4} \right)
\end{array}
\end{equation}
and
\begin{equation}
\begin{array}{rcl}
\phi_2(x) \!\!&\!\!=\!\!&\!\!
\left( 60\,{\pi }^{2}-128\,\pi +16 \right) {x}^{2}                                                            \\[1.25 ex]
\!\!&\!\!+\!\!&\!\!
\left( -65\,{\pi }^{3}+174\,{\pi }^{2}+12\,\pi -120 \right) {x}                                               \\[1.25 ex]
\!\!&\!\!+\!\!&\!\! \left(30\,{\pi }^{4}-122\,{\pi }^{3}-64\,{\pi
}^{2}+240\,\pi\right)
\end{array}
\end{equation}
are quadratic trinomials. Let us denote by $y_1$ the minimum of the trinomial $\phi_1(x)$
and by $y_2$ the minimum of the trinomial $\phi_2(x)$ over $[0, c_1]$ respectively.
It then becomes possible to verify that $y_1>271$ and $y_2 = \phi_2(c_1)>-815$. Thus
\begin{equation}
\label{ineq_P_6}
\begin{array}{rcl}
P(x)\!\! &\geq&\!\!y_1 + 4x^3\left(\pi-2\right)\left(\left(-20\pi+40 \right)x^3+y_2\right)                                                           \\[2ex]
     \!\!&>   &\!\!271 + 4x^3\left(\pi-2\right)\left(\left(40-20\pi \right)x^3-815\right)                                                            \\[2ex]
     \!\!&>   &\!\!271 - 4\left(\frac{23}{100}\right)^{\!3}\left(\pi-2\right)\left(\left(20\pi-40\right)\left(\frac{23}{100}\right)^{\!3}+815\right) \\[2ex]
     \!\!&>   &\!\!225 \,>\, 0,
\end{array}
\end{equation}
for $x \!\in\! \left( 0, \frac{23}{100} \right)$. Therefore $P(x)>0$
for $x \!\in\! \left( 0, c_1 \right) \!\subset\!
\left( 0, \frac{23}{100} \right)$.

\medskip
\noindent Now we prove that:
\begin{equation}
g_{2}(x) = P(x)\cos^2\!x - Q(x)\sin^2\!x > 0
\end{equation}
for $x \!\in\! (0,c_1)$. Let us note that $g_{2}(x)$ is a mixed trigonometric polynomial function,
and that the proof of previous inequality will be proved applying the methods from
\cite{Mortici_2011} and \cite{Malesevic_Makragic_2015}. In particular,
we use the following inequalities from \cite{Malesevic_Makragic_2015}:
\begin{equation}
\cos x>1-\frac{x^2}{2}, \qquad {\big (}x \in (0,c_1)\subseteq
(0,\sqrt{30}){\big )},
\end{equation}
and

\begin{equation}
\sin x<x, \qquad {\big (}x \in (0,c_1)\subseteq (0,\sqrt{20}){\big
)}.
\end{equation}
Therefore

\begin{equation}
g_{2}(x) > T_{10}(x) = P(x) \! \left( 1 - \frac{\;x^2}{2} \right)^2 \! - \, Q(x) \, x^2
\end{equation}

\smallskip\noindent
for $x \in (0, c_1)$ and it is enough to prove

\begin{equation}
\label{ineq_T_10} T_{10}(x) > 0,
\end{equation}

\smallskip \noindent
for $x \in (0, c_1)$. For the polynomial
\begin{equation}
\begin{array}{rcl}
T_{10}(x) \!\!&\!\!=\!\!&\!\!
\left( -20\,{\pi}^{2}\!+\!80\,\pi\!-\!80 \right) {x}^{10}                       \\[1.25 ex]
\!\!&\!\!+\!\!&\!\!
\left( 60\,{\pi}^{3}\!-\!248\,{\pi}^{2}\!+\!272\,\pi\!-\!32 \right) {x}^{9}     \\[1.25 ex]
\!\!&\!\!+\!\!&\!\! \left( -65\,{\pi}
^{4}\!+\!304\,{\pi}^{3}\!-\!336\,{\pi}^{2}\!-\!144\,\pi +240
\right) {x}^{8}                                                                 \\[1.25 ex]
\!\!&\!\!+\!\!&\!\!
 \left( 30\,{\pi}^{5}\!-\!182\,{\pi}^{4}\!+\!180\,{\pi}^{3}\!+\!400\,{\pi}^{2}
 -608\,\pi +128 \right) {x}^{7}                                                 \\[1.25 ex]
\!\!&\!\!+\!\!&\!\! \left(
-5\,{\pi}^{6}\!+\!53\,{\pi}^{5}\!-\!33\,{\pi}^{4}\!-\!392\,{\pi}^{3}
+424\,{\pi}^{2}\!+\!896\,\pi\!-\!1280 \right) {x}^{6}                           \\[1.25 ex]
\!\!&\!\!+\!\!&\!\! \left(
-6\,{\pi}^{6}\!-\!4\,{\pi}^{5}\!+\!190\,{\pi}^{4}\!+\!200\,{\pi}^{3}
-2464\,{\pi}^{2}\!+\!3008\,\pi\!-\!128 \right) {x}^{5}                          \\[1.25 ex]
\!\!&\!\!+\!\!&\!\! \left(
960+2400\,{\pi}^{3}\!-\!2784\,{\pi}^{2}\!-\!576\,\pi\!-\!413\,{\pi}^{4}
\!+\!11/4\,{\pi}^{6}\!-\!45\,{\pi}^{5} \right) {x}^{4}                          \\[1.25 ex]
\!\!&\!\!+\!\!&\!\! \left(
4\,{\pi}^{6}+196\,{\pi}^{5}-1136\,{\pi}^{4}+1200\,{\pi}^{3}
+1472\,{\pi}^{2}-1920\,\pi \right) {x}^{3}                                      \\[1.25 ex]
\!\!&\!\!+\!\!&\!\! \left(
-36\,{\pi}^{6}+264\,{\pi}^{5}-192\,{\pi}^{4}-1184\,{\pi}^{3}
+1440\,{\pi}^{2} \right) {x}^{2}                                                \\[1.25 ex]
\!\!&\!\!+\!\!&\!\!
\left( -24\,{\pi}^{6}-16\,{\pi}^{5}+408\,{\pi}^{4}-480\,{\pi}^{3} \right) x     \\[1.25 ex]
\!\!&\!\!+\!\!&\!\! 11\,{\pi}^{6}-52\,{\pi}^{5}+60\,{\pi}^{4}
\end{array}
\end{equation}
we form the polynomials:
\begin{equation}
\begin{array}{rcl}
\psi_{1}(x) \!&\!\!=\!\!&\! {\Big (} (
-20\,{\pi}^{2}\!+\!80\,\pi\!-\!80 ) x +
( 60\,{\pi}^{3}\!-\!248\,{\pi}^{2}\!+\!272\,\pi\!-\!32 ) {\Big )} x^{9},                          \\[1.5 ex]
\psi_{2}(x) \!&\!\!=\!\!&\!
{\Big (} ( -65\,{\pi} ^{4}\!+\!304\,{\pi}^{3}\!-\!336\,{\pi}^{2}\!-\!144\,\pi +240 ) x^2          \\[0.75 ex]
\!&\!\! \!\!&\! +
( 30\,{\pi}^{5}\!-\!182\,{\pi}^{4}\!+\!180\,{\pi}^{3}\!+\!400\,{\pi}^{2}-608\,\pi +128 ) x        \\[0.75 ex]
\!&\!\! \!\!&\! + (
-5\,{\pi}^{6}\!+\!53\,{\pi}^{5}\!-\!33\,{\pi}^{4}\!-\!392\,{\pi}^{3}+424\,{\pi}^{2}\!+\!896\,\pi\!-\!1280
)
{\Big )} {x}^{6},                                                                                 \\[1.5 ex]
\psi_{3}(x) \!&\!\!=\!\!&\! {\Big (} (
-6\,{\pi}^{6}\!-\!4\,{\pi}^{5}\!+\!190\,{\pi}^{4}\!+\!200\,{\pi}^{3}-2464\,{\pi}^{2}
\!+\!3008\,\pi\!-\!128) x^2                                                                       \\[0.75 ex]
\!&\!\! \!\!&\! + (
960+2400\,{\pi}^{3}\!-\!2784\,{\pi}^{2}\!-\!576\,\pi\!-\!413\,{\pi}^{4}\!+\!\frac{11}{4}\,{\pi}^{6}
\!-\!45\,{\pi}^{5} ) x                                                                           \\[0.75 ex]
\!&\!\! \!\!&\! +
(4\,{\pi}^{6}+196\,{\pi}^{5}-1136\,{\pi}^{4}+1200\,{\pi}^{3}
+1472\,{\pi}^{2}-1920\,\pi)
{\Big )} {x}^{3},                                                                                 \\[1.5 ex]
\psi_{4}(x) \!&\!\!=\!\!&\! {\Big (}
( -36\,{\pi}^{6}+264\,{\pi}^{5}-192\,{\pi}^{4}-1184\,{\pi}^{3}+1440\,{\pi}^{2} ) {x}^{2}          \\[0.75 ex]
\!&\!\! \!\!&\! +
( -24\,{\pi}^{6}-16\,{\pi}^{5}+408\,{\pi}^{4}-480\,{\pi}^{3} ) x                                  \\[0.75 ex]
\!&\!\! \!\!&\! + ( 11\,{\pi}^{6}-52\,{\pi}^{5}+60\,{\pi}^{4} )
{\Big )}.
\end{array}
\end{equation}

\vspace*{-2.5 mm}

\noindent
It is easy to check that
\begin{equation}
\hspace*{-4.5 mm}
\begin{array}{l}
\psi_1(x) > 0,                                                                                    \\[1.0 ex]
\psi_2(x) > 0,                                                                                    \\[1.0 ex]
\psi_3(x) > {\bigg (} (
-6\,{\pi}^{6}\!-\!4\,{\pi}^{5}\!+\!190\,{\pi}^{4}\!+\!200\,{\pi}^{3}-2464\,{\pi}^{2}
\!+\!3008\,\pi\!-\!128)\!\left(\frac{23}{100}\right)^2                                            \\[1.0 ex]
\;\;\;\;\;\;\;\;\;\;\;\;\;\; + (
960\!+\!2400\,{\pi}^{3}\!-\!2784\,{\pi}^{2}\!-\!576\,\pi\!-\!413\,{\pi}^{4}\!+\!\frac{11}{4}\,{\pi}^{6}
\!-\!45\,{\pi}^{5} )\!\left(\frac{23}{100}\right)                                                 \\[1.0 ex]
\;\;\;\;\;\;\;\;\;\;\;\;\;\; +
(4\,{\pi}^{6}\!+\!196\,{\pi}^{5}\!-\!1136\,{\pi}^{4}\!+\!1200\,{\pi}^{3}
\!+\!1472\,{\pi}^{2}-1920\,\pi)
{\bigg )}\!\left(\frac{23}{100}\right)^{3}                                                        \\[1.0 ex]
\;\;\;\;\;\;\;\;\;\;\;\;\;\;
> -27,                                                                                            \\[1.0 ex]
\psi_4(x) > {\bigg (} (
-36\,{\pi}^{6}+264\,{\pi}^{5}-192\,{\pi}^{4}-1184\,{\pi}^{3}+1440\,{\pi}^{2}
)\!\left(\frac{23}{100}\right)^2                                                                  \\[1.0 ex]
\;\;\;\;\;\;\;\;\;\;\;\;\;\; +
( -24\,{\pi}^{6}-16\,{\pi}^{5}+408\,{\pi}^{4}-480\,{\pi}^{3} )\!\left(\frac{23}{100}\right)       \\[1.0 ex]
\;\;\;\;\;\;\;\;\;\;\;\;\;\; + (
11\,{\pi}^{6}-52\,{\pi}^{5}+60\,{\pi}^{4} )
{\bigg )}\!\left(\frac{23}{100}\right)^{3}                                                        \\[1.0 ex]
\;\;\;\;\;\;\;\;\;\;\;\;\;\;
> 54,
\end{array}
\end{equation}

\vspace*{-1.5 mm}

\noindent
for $x \in (0, c_1) \subset \left(0, \frac{23}{100}\right)$. Thus we may conclude that
\begin{equation}
T_{10}(x) = \psi_1(x) + \psi_2(x) + \psi_3(x) + \psi_4(x) > 27 > 0,
\end{equation}
for $x \in (0,c_1)$. Based on $G''_{1}(x) = g_{2}(x) / \!\left(Q(x)\cos^2\!x\right)
> T_{10}(x) / \!\left(Q(x)\cos^2\!x\right)$, the inequality
\begin{equation}
G''_{1}(x)
> 0,
\end{equation}
was proved for $x \in (0,c_1)$. Let us notice that
\begin{equation}
\begin{array}{rcl}
G'_{1}(0) \!\!&\!\!=\!\!&\!\!
\nfrac{1}{5}\ln\!\left(\nfrac{\pi}{2}\right)-\nfrac{2\pi\!-\!6}{\pi}
\,=\,
\nfrac{1}{5}\ln\!\left(1\!+\!\left(\nfrac{\pi}{2}\!-\!1\right)\right)-\nfrac{2\pi\!-\!6}{\pi}
                                                                                \\[1.00 ex]
\!\!&\!\!>\!\!&\!\! \nfrac{1}{5}\displaystyle\sum\limits_{k=1}^{4}{
\left( \nfrac{\left(\frac{\pi}{2}\!-\!1\right)^{2k-1}}{2k-1} -
\nfrac{\left(\frac{\pi}{2}\!-\!1\right)^{2k}}{2k} \right)} -
\nfrac{2\pi\!-\!6}{\pi}                                                         \\[2.50 ex]
\!\!&\!\!=\!\!&\!\!
\mbox{\small $
-{\nfrac {{\pi }^{8}}{10240}}
\!+\!
{\nfrac {{\pi}^{7}}{560}}
\!-\!
{\nfrac {7\,{ \pi }^{6}}{480}}
\!+\!
{\nfrac {7\,{\pi}^{5}}{100}}
\!-\!
{\nfrac {7\,{\pi }^{4}}{32}}
\!+\!
{\nfrac {7\,{\pi }^{3}}{15}}
\!-\!
{\nfrac {7\,{\pi }^{2}}{10}}+{\nfrac{4\pi}{5}}
\!-\!
{\nfrac {3561}{1400}}
\!+\!
{\nfrac{6}{\pi}}
$} > 0
\end{array}
\end{equation}
and
\begin{equation}
G_{1}(0) = 0.
\end{equation}

\break

\noindent
Based on Theorem \ref{TH} it follows that $G_{1}(x) > 0$ for $x \in (0,c_1)$.
This also proves that $G(x) > 0$ for $x \in (0,c_1)$, which in turn proves that
$F(x) > 0$ for $x \in \left(c,\frac{\pi}{2}\right)$.

\smallskip
\noindent
Therefore, we can conclude that $F(x)>0$ for any $x \in \left(0, \frac{\pi}{2}\right)$.
The proof of Nishizava's open problem is now completed.

\section{Conclusions}

This paper proved an open problem stated by Nishizawa in \cite{Nishizawa_2015},
applying computation method from \cite{Mortici_2011} and \cite{Malesevic_Makragic_2015}.
We note that proofs of polynomial inequalities (\ref{A}), (\ref{P14}), (\ref{ineq_P_6})
and (\ref{ineq_T_10}) can be based on reducing (by differentiation) of the corresponding
polynomials to polynomials of a degree up to four (as~illustrated in papers
\cite{Banjac_Makragic_Malesevic_2015}\mbox{\small $-$}\cite{Malesevic_Banjac_Jovovic_2015}),
which allows symbolic radical representation of roots.

\smallskip
Our approach, based on the fact
(\ref{Properties}), allows new proofs of some power-exponential inequalities from papers
\cite{Nishizawa_2015}, \cite{Chen_Cheung_2011}\mbox{\small $-$}\cite{Banjac_Makragic_Malesevic_2015}
and monographs \cite{Mitrinovic_1970}, \cite{Milovanovic_Rassias_2014}.

\medskip

\bigskip
\noindent
\textbf{Acknowledgement.}
The authors are grateful to Professor G. Milovanovi\' c and Professor C. Mortici
for valuable comments on earlier draft of the paper.

\break

\end{document}